\documentclass[10pt]{article}
\usepackage{amssymb,amsmath,amsthm,latexsym}
\usepackage{amsfonts}
\usepackage[twlrm]{rawfonts}
\usepackage[english]{babel}
\usepackage[latin1]{inputenc}
\usepackage{fancyhdr}
\usepackage{indentfirst}
\usepackage{color}
\newtheorem{defi}{Definition}[section]

\newtheorem{teo}{Theorem}[section]
\newtheorem{cor}{Corolary}[section]
\newtheorem{lem}{Lemma}[section]
\newtheorem{pro}{Proposition}[section]
\newtheorem{claim}{Claim}[section]
\newenvironment{dem}[1][Proof]{\noindent\textbf{#1.} }{\hfill \rule{0.5em}{0.5em}}

\newcommand{\N}{\mathbb{N}}

\newcommand{\R}{\mathbb{R}}

\newcommand{\C}{C(\overline{\Omega})}

\newtheorem{ob}{Remark}[section]

\begin{document}
	
\setlength{\baselineskip}{6.5mm} \setlength{\oddsidemargin}{8mm}
\setlength{\topmargin}{-3mm}

\title{\bf Existence of solution for a nonlocal dispersal model with nonlocal term via bifurcation theory}

\author { Claudianor O. Alves$^a$\thanks{C. O. Alves was partially supported by CNPq/Brazil
 304804/2017-7, coalves@mat.ufcg.edu.br}\, , \, \ Natan de Assis Lima$^b$\thanks{N. A. Lima, natan.mat@cche.uepb.edu.br}\,\,\, and  \,    Marco A. S. Souto$^a$\thanks{M. A. S. Souto was partially supported by CNPq/Brazil 306082/2017-9 and INCT-MAT, marco@mat.ufcg.edu.br}\,\,\,\,\,\,\vspace{2mm}
\and {\small $a.$ Universidade Federal de Campina Grande} \\ {\small Unidade Acad\^emica de Matem\'{a}tica} \\ {\small CEP: 58429-900, Campina Grande - Pb, Brazil}\\
{\small $b.$ Universidade Estadual da Paraíba} \\ {\small Centro de Ciências Humanas e Exatas} \\ {\small CEP: 58500-000, Monteiro - Pb, Brazil}}


\date{}

\maketitle

\begin{abstract}
	In this paper we study the existence of solution for the following class of nonlocal problems
\[
L_0u =u \left(\lambda - \int_{\Omega}Q(x,y) |u(y)|^p dy \right) , \ \mbox{in} \ \Omega,
\]
where $\Omega \subset \R^{N}$, $N\geq 1$, is a smooth bounded domain, $p>0$, $\lambda$ is a real parameter, $Q:\Omega \times \Omega \to \R$ is a nonnegative function, and $L_0 : \C \to \C$ is a nonlocal dispersal operator. The existence of solution is obtained via bifurcation theory.


\end{abstract}

\section{Introduction}

In this work we study the existence of positive solution for the following equation 
\[
L_0u =u \left(\lambda - \int_{\Omega}Q(x,y) |u(y)|^p dy \right) , \ \mbox{in} \ \Omega, \eqno{(P)}
\]
where $\Omega \subset \R^{N}$, $N\geq 1$, is a smooth bounded domain, $p>0$, $\lambda$ is a real parameter, $Q:\Omega \times \Omega \to \R$ is a nonnegative function with $Q \in C(\overline{\Omega} \times \overline{\Omega})$ and verifying some hypotheses that will be detailed below, and $L_0 : \C \to \C$ is the nonlocal dispersal operator given by 
\begin{equation}\label{eq02}
L_0 u(x) = \int_{\Omega}K(x,y)u(y)dy, \quad \mbox{for} \ u \in \C, \end{equation}
with a continuous and nonnegative dispersal kernel $K$. The dispersion mechanism is currently a focus of theoretical interest and has received much attention recently. Most of these continuous dispersion models are based on reaction-diffusion equations, which are widely studied see \cite{Bates-Chmaj}, \cite{Bates-Fife-Ren-Wang}, \cite{Bates-Zhao}, \cite{Chasseigne-Chaves-Rossi}, \cite{Chen}, \cite{Coville 2}, \cite{Coville 3}, \cite{Garcia-Rossi}, \cite{Garcia-Rossi 2}, \cite{Kao-Lou-Shen}, \cite{Kao-Lou-Shen 2}. This type of diffusion process has been widely used to describe the dispersion of a population (of cell or organisms) through the environment, as indicated in \cite{Fife 1}, \cite{Fife 2}, \cite{Hutson-Martinez-Mischaikow-Vickers}, if $u(y)$ is thought of as a density at a location $y$, $K(x,y)$ as the probability distribution of jumping from a location $y$ to a location $x$, then the rate at which the individuals from all other places are arriving at location $x$ is $$\int_{\Omega}K(x,y)u(y)dy.$$  In this context, $\lambda$ is a parameter which represents the intrinsic growth rate of the species, and the nonlocal term $$\int_{\Omega}Q(x,y)|u(y)|^pdy$$ can be interpreted as a weighted average of $u$ at
all the domain. In many problems in biology (and ecology), for example seed dispersal problems, this formulation of the dispersion of individuals finds its justification; see \cite{Cain-Milligan-Strand}, \cite{Clark}, \cite{Medlock-Kot}, \cite{Murray}, \cite{Schurr-Steinitz-Nathan}. The presence of nonlocal reaction term in equation $(P)$ means, from the biological point of view, that the crowding effect depends not only on their own point in space but also depends on the entire population in an $N$-dimensional habitat $\Omega$, see \cite{Furter-Grinfeld}.

The motivation to study $(P)$ comes from the model to study the behavior of a species inhabiting in a smooth bounded domain $\Omega$, whose the classical logistic equation with laplacian diffusion is given by
\begin{equation}\label{eq03}
\left\{
\begin{array}{lcl}
-\Delta u = u \left(\lambda - b(x)u^p \right)\quad \mbox{in} \ \Omega, \\
u=0 \quad \mbox{on} \ \partial \Omega,
\end{array}
\right.
\end{equation}
where $b(x)$ describes the limiting effect of crowding of the population. In (\ref{eq03}), we are assuming that $\Omega$ is surrounded by inhospitable areas, due to the homogeneous Dirichlet boundary conditions. Note that the equation in (\ref{eq03}) is a local equation, and so the crowding effect of the
population $u$ at $x$ only depends on the value of the population in the same point $x$. In \cite{Chipot}, Chipot has considered that crowding effect depends also on the value of the population around of $x$, that is, the crowding effect depends on the value of $u$ in a neighborhood of $x$, $B_r(x)$, the centered ball at $x$ of radius $r>0$. To be more precisely, Chipot considered the nonlocal problem 
\begin{equation}\label{eq04}
\left\{
\begin{array}{lcl}
-\Delta u = u \left(\lambda - \int_{\Omega \cap B_r(x)}b(y)u^p(y)dy \right)\quad \mbox{in} \ \Omega, \\
u=0 \quad \mbox{on} \ \partial \Omega,
\end{array}
\right.
\end{equation}
where $b$ is a nonnegative and nontrivial continuous function. After that, a special attention was given to the problem 
\begin{equation}\label{eq05}
\left\{
\begin{array}{lcl}
-\Delta u = u \left(\lambda - \int_{\Omega}Q(x,y)u^p(y)dy \right)\quad \mbox{in} \ \Omega, \\
u=0, \quad \mbox{on} \ \partial \Omega, 
\end{array}
\right.
\end{equation}
by supposing different conditions on $Q$, see for example \cite{Allegretto-Nistri}, \cite{Alves-Delgado-Souto-Suarez}, \cite{Chen-Shi}, \cite{Correa-Delgado-Suarez}, \cite{Coville}, \cite{Leman-Meleard-Mirrahimi} and \cite{Sun-Shi-Wang} and their references.

In \cite{Chen-Shi}, Chen and J. Shi have considered the case $p=1$ and the kernel function $Q(x,y)$ being a continuous and nonnegative function on $\Omega \times \Omega$ with  $\displaystyle\int_{\Omega}Q(x,y)u(y)dy >0$ for all positive continuous functions $u$ on $\Omega$. In that paper was proved the existence of $\lambda^* > \lambda_1$ such that (\ref{eq05}) possesses at least a positive solution for $\lambda \in (\lambda_1, \lambda^*]$.

In \cite{Allegretto-Nistri}, Allegretto and P. Nistri have showed that (\ref{eq05}) possesses a unique positive solution when $\lambda> \lambda_1$ and $Q(x,y)= Q_\delta (|x-y|)$ is a mollifier in $\R^N$, i.e., $Q_\delta (|x-y|) \in C^{\infty}_{0}(\overline{\Omega})$, $\displaystyle\int_{\R^N}Q_\delta (|x-y|)dy =1$ for any $x$ with $$Q_\delta (|x-y|) =0 \quad |x-y|\geq \delta$$ and $$Q_\delta (|x-y|) \quad \mbox{bounded away from zero is} \quad |x-y|<\mu < \delta.$$ Observe that in this case, $Q$ vanishes away from the diagonal of $\Omega \times \Omega$.

In \cite{Sun-Shi-Wang}, Sun, Shi and Wang have investigated the existence of positive solutions for  (\ref{eq05}) with $Q(x,y)= Q_1(|x-y|)$ and $\Omega = (-1,1)$, where $Q_1 : [0,2] \to (0,\infty)$ is a nondecreasing and piecewise continuous function satisfying $$\int_{0}^{2}Q_1(y)dy > 0.$$

When $Q(x, y)$ is a separable variable, i.e.,
$Q(x, y) = g(x)h(y)$ with $h \geq 0$; $h \neq 0$ and $g(x) > 0$, Corrêa, Delgado and Suárez \cite{Correa-Delgado-Suarez} have studied (\ref{eq05}) and proved the existence and uniqueness of positive solution. Moreover, in  Coville \cite{Coville} and  Leman, Méléard and Mirrahimi \cite{Leman-Meleard-Mirrahimi}, by assuming $g \equiv 1$, $p > 1$ and homogeneous Neumann boundary conditions, the authors have proved that positive solution of (\ref{eq05}) attracts all the possible solutions of the corresponding parabolic associated with (\ref{eq05}). When $g \geq 0$, $g \neq 0$ and $g \equiv 0$ in $\Omega_0 \subset \Omega$, then (\ref{eq05}) possesses a unique positive solution for $\lambda \in (\lambda_1, \lambda_0)$ where $\lambda_0$ is the principal eigenvalue of the minus laplacian in $\Omega_0$, for more details see \cite{Correa-Delgado-Suarez}.

Finally, in \cite{Alves-Delgado-Souto-Suarez}, Alves, Delgado,  Souto and Su\'{a}rez have considered the existence and nonexistence of solution for (\ref{eq05}). In that paper, they have studied a more general problem than the previous ones, more precisely, they have considered that $Q$ satisfies:
\begin{itemize}
\item [$(Q_{1})$] $Q \in L^{\infty}(\Omega \times \Omega)$ and $Q(x,y) \geq 0$ for all $x,y\in \Omega$;
\item [$(Q^{'}_2)$] If $w$ is measurable and $\displaystyle\int_{\Omega \times \Omega}Q(x,y)|w(y)|^{p}|w(x)|^2dxdy =0$, then $w =0$ a.e. in $\Omega$.
\end{itemize}

In the present paper, our main goal is showing the existence of solution for $(P)$ via the global bifurcation theorem result due to Rabinowitz, see Theorem $29.1$ in \cite{Deimling}. The difference of the problem with $L_0$ of the problem with $-\Delta$ is that, $(-\Delta)^{-1}$ is a compact operator, while $(L_0+M)^{-1}$ is not a compact for any  $M \geq 0$ large. This difference between the operators brings some difficulties and the bifurcation here is made in another way.

This paper is organized as follows. In Section 2 we show a Krein-Rutman type result for $L_0$. In Section 3 we study the existence of a positive solution for equation $(P)$ without boundary conditions, by supposing that $K$ satisfies:
\begin{itemize}
\item [$(K_{1})$] $K(x,y)=K(y,x)$ for all $x,y\in \overline\Omega$;
\item [$(K_{2})$] There exists $\delta>0$ such that $K(x,y) > 0$ for all $x,y \in \overline{\Omega}$ and $|x-y|\leq \delta$;
\end{itemize}
and supposing that $Q$ is a continuous function that satisfies:
\begin{itemize}
\item[$(Q_2)$] There exist $r,\sigma>0$ such that $Q(x,y) \geq \sigma$ for all $x,y \in \overline{\Omega}$ and $|x-y|\leq r$.
\end{itemize}
If the inequality above works for all $x,y$ in $\Omega$, we say that $Q$ satisfies $(Q''_2)$, that is, 
\begin{itemize}
	\item[$(Q''_2)$] $Q(x,y) \geq \sigma$ for all $x,y \in \overline{\Omega}$.
\end{itemize}

The conditions $(K_1)$ and $(K_2)$ above are hypotheses usually considered on the operator $L_0$, as we can see in \cite{Bates-Zhao}, \cite{Garcia-Rossi} and \cite{Garcia-Rossi 2}.

We would like point out that assumption $(Q_2)$ implies in $(Q^{'}_2)$.

\begin{defi}
For each $Q: \overline{\Omega} \times \overline{\Omega} \longrightarrow \R$ we define the oscillation of $Q$ in $x$, uniformly in $y$, by
$$[Q] = \sup_{x,y,z \in \Omega}|Q(x,y) -Q(z,y)|.$$
\end{defi}

With the above hypotheses we have our first result that establishes the existence of local bifurcation. 

\begin{teo} \label{Teo 1}
Suppose that  $p>0$, $[Q]>0$, $(K_1)-(K_2)$ and $(Q''_2)$ hold.
Then the problem $(P)$ has a positive solution for all $\lambda \in (\lambda_1, \lambda_1 + \displaystyle\frac{\lambda_1 \sigma}{[Q]})$, where $\lambda_1$ is the principal eigenvalue of  $L_0$.
\end{teo}

It is a corollary of the proof of the Theorem \ref{Teo 1} that if $[Q]=0$ we have solution for all $\lambda >\lambda_1$.

In order to obtain a global bifurcation result we will assume the following assumption on $Q$: 
\begin{itemize}
\item[$(Q_3)$]  There are $x_0 \in \overline{\Omega}$ and a nonnegative function $a:\overline{\Omega} \longrightarrow\R$, such that $a^{-1} \in L^{q}(\Omega)$ where $q=\max\{1,p\}$ and $Q(x_0,y) \geq Q(x,y) + a(x)$ for all $x,y \in \Omega$.
\end{itemize}

Note that $(Q_3)$ implies that $Q(x_0,y)\geq Q(x,y)$, for all $x,y \in \Omega$ and $a(x_0)=0$. Here is an example of a function $Q$ that satisfies $(Q_3)$:
$$Q(x,y)=h(y)[M-|x-x_1|^{q_1} |x-x_2|^{q_2} ... |x-x_k|^{q_k}] +g(y)$$ where $x_i \in \overline{\Omega}$, $q_i<\displaystyle\dfrac{N}{p}$, $M>0$ is large enough, $h$ and $g$ are positive continuous functions on $\overline{\Omega}$. It is easy to see that $(Q_3)$ works for $x_0$ being any $x_i$ and $a(x)=m|x-x_1|^{q_1} |x-x_2|^{q_2} ... |x-x_k|^{q_k}$, for $m>0$ small enough.

The next result guarantees the existence of a connected component of solutions which contains solution for our problem $(P)$ for any $\lambda>\lambda_1$.

\begin{teo} \label{Teo 2}
Assume $p>0$, $(K_1)-(K_2)$ and $(Q_2)-(Q_3)$ hold. Then, there is a connected component of positive solutions of $(P)$ coming out from $\lambda_1>0$, with the property that includes solutions of the form $(\lambda, u)$ for all $\lambda > \lambda_1$.
\end{teo}

In Section 4, under a weaker condition on $Q$, we can solve the problem for all $\lambda>\lambda_1$, but we can not guarantee the existence of a connected component of positives solutions for $(P)$. Here we consider that $Q$ satisfies:
\begin{itemize}
\item[$(Q_4)$]  There are $x_0 \in \overline{\Omega}$ such that $Q(x_0,y) \geq Q(x,y)$ for all $x,y \in \Omega$.
\end{itemize}

The main result this paper is the following:

\begin{teo} \label{Teo 3}
Assume $p>0$, $(K_1),(K_2)$, $(Q_2)$ and $(Q_4)$ hold. Then, problem $(P)$ has a positive solution for all $\lambda > \lambda_1$. 
\end{teo}

In fact, we comment in the final of the  Section 4 that we have the same result of Theorem \ref{Teo 3} if $Q$ satisfies a more general condition than $(Q_4)$: \begin{itemize}
	\item[$(Q'_4)$]  There are a decomposition $\Omega=\displaystyle{\bigcup_{j=1}^m E_j}$, $x_1\in E_1$, $x_2\in E_2$,... $x_m\in E_m$  such that $Q(x_j,y) \geq Q(x,y)$ for all $x\in E_j$, $y \in \Omega$.
\end{itemize}

\section{On the principal eigenvalue of $L_0$}

In this section, we consider some preliminary facts related to the principal eigenvalue of $L_0$. Let us denote $\C$ by $X$ and consider the dispersal operator $L_0: X \to X$, where kernel $K$ is a nonnegative continuous function verifying $(K_1)$. It is easy to see that $L_0$ is a compact operator and $L_0(X_+) \subset X_+$, where $X_+$ is the positive cone in $X$, that is, $$X_+=\{u \in X; u(x)\geq 0, \forall x \in \overline{\Omega}\}.$$ Furthermore, we can also consider $L_0:L^{2}(\Omega) \to L^{2}(\Omega)$, which is well defined, compact and symmetric with $L_0(L^{2}(\Omega)) \subset X$. It is well known that the resolvent of $L_0$ is defined by $$\rho(L_0)=\{\lambda \in \R; L_0-\lambda I \ \mbox{is bijective}\}$$ and its spectrum is $\sigma(L_0)=\R \backslash \rho(L_0)$. By spectral theory of compact operators, we have for $\lambda \neq 0$
 \[\lambda \in \sigma(L_0) \iff N(L_0-\lambda I)\neq \{0\} \mbox{ or } R(L_0-\lambda I)= (L_0-\lambda I)(X)\neq X.  \]
 Here, $EV(L_0)$ denotes the eigenvalues set of $L_0$ given by
  \[
  EV(L_0)=\{\lambda \in \mathbb R: N(L_0-\lambda I)\neq \{0\} \}.
  \]
 To avoid some confusion let us denote by $\tilde\sigma(L_0)$ the spectrum of $L_0:L^2(\Omega)\to L^2(\Omega)$, and $\tilde{EV}(L_0)$ its eigenvalues set. Note that, $\sigma(L_0)=\tilde\sigma(L_0)$, that is $$\lambda\in \tilde{EV}(L_0) \iff\lambda\in EV(L_0).$$ Indeed, from $X\subset L^2(\Omega)$ we have that $EV(L_0)\subset \tilde{EV}(L_0)$,  and the sufficient condition is done. Now suppose that $\lambda$ is an eigenvalue of $L_0:L^2(\Omega)\to L^2(\Omega)$, then there is $w\in  L^2(\Omega)\setminus\{0\}$ such that 
 \[
 \lambda w=L_0w \in L_0(L^2(\Omega))\subset X \Rightarrow \tilde{EV}(L_0)\subset {EV}(L_0).
 \]
 
 Since $L_0$ is a symmetric operator in $L^2(\Omega)$, that is, 
 \[
 \langle L_0u,v \rangle= \langle u,L_0v\rangle,\quad \forall u,v\in L^2(\Omega),
 \]
 where $ \langle u,v\rangle=\int_\Omega uv dx$ is the inner product of $L^2(\Omega)$, we have $\tilde{\sigma}(L_0)\subset [m_o,m]$, where
 \[
 m_o=\inf_{u\in L^2(\Omega)\setminus\{0\}}\frac{\langle L_0u,u\rangle}{\int_\Omega |u|^2dx} \mbox { and } m=\sup_{u\in L^2(\Omega)\setminus\{0\}}\frac{\langle L_0u,u\rangle}{\int_\Omega |u|^2dx}.
 \]
 Moreover $m_o,m\in \tilde{\sigma}(L_0)$ (see Brézis \cite[Proposition 6.9, pg 165]{Brezis}), and so, $ m=\sup \tilde{\sigma}(L_0).$ From definition $m$ and positiveness of $K$, it follows that 
 \[
 0<m=\sup_{u\in L^2(\Omega)\setminus\{0\}}\frac{\int_{\Omega\times \Omega}K(x,y)u(x)u(y)dxdy} {\int_\Omega |u|^2dx}
 \]
and there is $w \in L^{2}(\Omega) \setminus \{0\}$ such that 
\[
m=\frac{\int_{\Omega\times \Omega}K(x,y)w(x)w(y)dxdy}{\int_\Omega |w|^2dx}.
\]	 
Thus, $m \in \tilde{EV}(L_0)$ and $w$ is an eigenfunction of $L_0$ associated with $m$.

 The next result establishes that if $w \in L^{2}(\Omega)$ is an eigenfunction  associated with the eigenvalue $m$, then $w$ does not change sign.
 
 \begin{lem} \label{1}
 	Suppose that $w\in L^2(\Omega)$ is such that
 	\[ m= \frac{\int_{\Omega\times \Omega}K(x,y)w(x)w(y)dxdy}{\int_\Omega |w|^2dx}.\]
 	Then $w$ is continuous and $L_0w=mw$, that is, $m$ is the maximum eigenvalue of $L_0$. Moreover, since $\Omega$ is a connected set, we must have $w>0$ in $\Omega$ or $w<0$ in $\Omega$ ($w$ does not change signal).
 \end{lem}

\begin{dem}
To begin with, we show that $L_0w=mw$. In fact, by definition of $m$, for any $t\in \mathbb R$ and $v\in L^2(\Omega)$, we must have
\[
\langle L_0(w+tv),w+tv \rangle \leq m\int_\Omega (w+tv)^2dx
\]
that is,
\[\langle L_0w,w \rangle +2t\langle L_0w,v \rangle+t^2\langle L_0v,v \rangle \leq m\int_\Omega w^2dx+2mt\int_\Omega wvdx+mt^2\int_\Omega v^2dx,\]
\[2t\langle L_0w,v \rangle+t^2\langle L_0v,v \rangle \leq 2mt\int_\Omega wvdx+mt^2\int_\Omega v^2dx,\]
which yields 
\[
\langle L_0w,v \rangle+\frac t2\langle L_0v,v \rangle \leq m\int_\Omega wvdx+\frac {mt}2\int_\Omega v^2dx, \mbox { if } t>0
\]
and
\[
\langle L_0w,v \rangle+\frac t2\langle L_0v,v \rangle \geq m\int_\Omega wvdx+\frac {mt}2\int_\Omega v^2dx, \mbox { if } t<0.
\]
Taking the limit as $t\to 0$, we get $\langle L_0w,v \rangle=m\int_\Omega wvdx$, for all $v\in L^2(\Omega)$, from where it follows that, $L_0w=mw$. Since $L_0w$ is a continuous function, the equality $L_0w=mw$ implies that $w$ is a continuous function on $\overline{\Omega}$.

In the sequel, we prove that if $w\neq0$ and $w\geq 0$ in $\Omega$, then $w>0$ in $\Omega$. Indeed, fixing $C=w^{-1}(\{0\})=\{x\in \Omega:w(x)=0\}$, we have that $C$ is an open set in $\Omega$, because if $z \in C$,  we have that 
$$
0 \leq \int_{\Omega \cap B_{\delta}(z)}K(z,y)w(y)dy \leq \int_{\Omega}K(z,y)w(y)dy = m w(z) = 0,
$$ 
implying that $w(x)=0$ for all $ x \in \Omega \cap B_{\delta}(z)$, that is, $ \Omega \cap B_{\delta}(z) \subset C$. This proves that $C$ is an open set in $\Omega$. Moreover, from continuity of $w$, it is easy to see that $C$ is also a closed set in $\Omega$. Recalling that $\Omega$ is connected and $\Omega \setminus C$ is non-empty, we conclude that $w>0$ in $\Omega$.

To finish the proof,  we will show that $w$ does not change signal. To this end, we fix the sets
\[A=\{x\in \Omega: w(x)>0\} \mbox { and }B=\{x\in \Omega: w(x)<0\}. \]
If $w$ changes signal, both subsets $A$ and $B$ have positive Lebesgue measure. By a simple computation, we see that
\begin{eqnarray*}
\int_{\Omega\times \Omega}K(x,y)w(x)w(y)dxdy=\int_{A\times A}K(x,y)w(x)w(y)dxdy \\ +\int_{B\times B}K(x,y)w(x)w(y)dxdy\\
+\int_{A\times B}K(x,y)w(x)w(y)dxdy\\
+\int_{B\times A}K(x,y)w(x)w(y)dxdy.
\end{eqnarray*}
If $K(x,y)$ is not zero in $A\times B$, we derive that 
\begin{eqnarray*}
	\int_{A\times A}K(x,y)w(x)w(y)dxdy=\int_{A\times A}K(x,y)|w(x)|\,|w(y)|dxdy,\\
	\int_{B\times B}K(x,y)w(x)w(y)dxdy=\int_{B\times B}K(x,y)|w(x)|\,|w(y)|dxdy\\ 
	\int_{A\times B}K(x,y)w(x)w(y)dxdy<0<\int_{A\times B}K(x,y)|w(x)|\,|w(y)|dxdy
\end{eqnarray*}
and 
\begin{eqnarray*}
	\int_{B\times A}K(x,y)w(x)w(y)dxdy<0<\int_{B\times A}K(x,y)|w(x)|\,|w(y)|dxdy.
\end{eqnarray*}
These informations lead to 
\[
 m= \frac{\int_{\Omega\times \Omega}K(x,y)w(x)w(y)dxdy}{\int_\Omega w^2dx}<  \frac{\int_{\Omega\times \Omega}K(x,y)|w(x)|\,|w(y)|dxdy}{\int_\Omega |w|^2dx}\leq m,
\]
which is impossible, then $w> 0 $ in $\Omega$ or $w <0 $ in $\Omega$.

Suppose now that $K(x,y)\equiv 0$ in $(A\times B)\cup(B\times A)$. Then, 
\begin{eqnarray*}
	a_1=\int_{A\times A}K(x,y)w(x)w(y)dxdy,\\
		a_2=\int_{B\times B}K(x,y)w(x)w(y)dxdy,
\end{eqnarray*}
and
\begin{eqnarray*}		
		b_1=\int_A  w^2dx \mbox{ and } 	b_2=\int_B w^2dx.
\end{eqnarray*}
Consequently
\[
m=\frac{a_1+a_2}{b_1+b_2},\, \,\frac{a_1}{b_1}\leq m \mbox { and } \frac{a_2}{b_2}\leq m.
\]
A simply computation shows that $a_1=mb_1$, and so, the function $w^+=\max\{w,0\}$ satisfies 
\[ 
m= \frac{\int_{\Omega\times \Omega}K(x,y)w^+(x)w^+(y)dxdy}{\int_\Omega (w^+)^2dx}.
\] 
By using the first part of this proof, we must have $w^+>0$ in $\Omega$, or equivalently, $w(x)>0$ in ${\Omega}$, contradicting the fact that $|B|>0$, finishing the proof. 
\end{dem}

\begin{cor}
	If $w$ is an eigenfunction of $L_0$ associated to $m$ then $w$ must be a positive (or negative) eigenfunction.  Besides, $\dim N(L_0-mI)=1$.
\end{cor}

\begin{dem}
If $L_0w=mw$, then $\langle L_0w,w\rangle=m \displaystyle\int_\Omega w^2$ and $w\neq0$. From Lemma \ref{1}, $w$ must be a positive (or negative) eigenfunction. For the second part of the corollary, suppose that there are two eigenfunctions $w,\phi$ associated with $m$ that are linearly independent. Then, without loss of generality we can suppose that $\displaystyle\int_\Omega w\phi dx=0$. However, it is impossible, because $w$ and $\phi$ has defined signal in $\Omega$, then $\displaystyle\int_\Omega w\phi dx \not= 0$.
\end{dem}
 
\begin{cor} \label{2}
Under the conditions $(K_1)-(K_2)$,  if $w$ is an eigenfunction of $L_0$ associated to $m$, then $w$ is positive in $\overline{\Omega}$ (or negative in $\overline{\Omega}$). Hence, $w$  is discontinuous  on $\partial \Omega$ if we define $w=0$ in $\mathbb R^N\setminus \Omega$. 
\end{cor}
\begin{dem}
Indeed, if $x_0 \in \partial \Omega$ is such that $w(x_0)=0$, then 
$$
0<\int_{\Omega \cap |x_0-y|\leq \delta}K(x_0,y)w(y)dy \leq \int_{\Omega}K(x_0,y)w(y)dy = \lambda_1 w(x_0)=0
$$ 
which is absurd, showing the desired result.  
\end{dem}

As a byproduct of the study made until moment, we have the following result

\begin{pro} \label{3}
	The eigenvalue problem 
	\[
	L_0 u = \lambda u, \quad \mbox{in} \ \Omega, 
	\]
	has an unique eigenvalue $\lambda_1 > 0$ whose the eigenfunction are continuous on $\overline{\Omega}$ with defined signal and $\dim N(L-\lambda_1 I)=1$. Moreover, $\lambda_1 =m=\sup \sigma(L_0)$.
\end{pro}

Now, before proving our next result, it is necessary to recall that assumptions on $K$ implies that for each $u \in \C$, with $u \geq 0$ in $\Omega$ only one of the possibilities below holds:
$$L_0u >0 \ \mbox{in} \ \overline{\Omega} \ \mbox{or} \ u\equiv 0 \ \mbox{in} \ \overline{\Omega}.$$ Hence, we have the following lemma

\begin{lem} \label{4}
Suppose that $u \in \C$, $u \geq 0$, $u \neq 0$ in $\Omega$ and $c(x)$ given by $L_0u=c(x)u$, then $\|c\|_\infty \geq \lambda_1$. This inequality becomes equality only when $c(x) \equiv \lambda_1$.
\end{lem}

\begin{dem}
Consider $\varphi_1$ a positive eigenfunction associated to $\lambda_1$. Therefore, $L_0 \varphi_1 = \lambda_1 \varphi_1$ and 
\[
\lambda_1 \int_{\Omega}u\varphi_1 dx = \int_{\Omega}uL_0\varphi_1 dx = \int_{\Omega}\varphi_1L_0u dx = \int_{\Omega}c(x)u\varphi_1 dx \leq \|c\|_\infty \int_{\Omega}u\varphi_1 dx.
\]
Since $\displaystyle\int_{\Omega}u\varphi_1 dx >0$ the above inequality leads to $\|c\|_\infty \geq \lambda_1$. 
\end{dem}

\begin{lem} \label{5}
Let $u \in L^1(\Omega)$ be a nonnegative function and $c \in \C$ satisfying $$L_0u(x)=c(x)u(x) \quad \mbox{a.e. in} \ \Omega.$$ Then, if $u \neq 0$ we have that $u$ is continuous and positive in $\overline{\Omega}$. Furthermore, $c$ is positive in $\overline{\Omega}$.
\end{lem}

\begin{dem}
Clearly $L_0u \in \C$,  and so,  $c(\cdot)u \in \C$. Consider the following sets, $$V=\{x \in \overline{\Omega};\ \mbox{there exists a ball}\ B\ \mbox{centered at point} \ x\ \mbox{such that}\ u \equiv 0\ \mbox{a.e. in}\ B \cap \Omega \}$$ and $W=\{x \in \overline{\Omega}; L_0u(x)>0\}$. Both subsets $V$ and $W$ are open in $\overline{\Omega}$. Now, we are showing that $W \cap V = \emptyset$ and $V=\overline{\Omega} \setminus W$. 

Indeed, if $z \notin W$ we have $L_0u(z)=0$. Thus, $$0=L_0u(z)=\int_{\Omega} K(z,y)u(y)dy\geq \int_{B_\delta(z)\cap \Omega}K(z,y)u(y)dy,$$ since $K(z,y)>0$ for all $|z-y|\leq \delta$, we get $u(y)=0$ a.e. in $B_\delta (z)\cap \Omega$, that is $z \in V$. Since $\Omega$ is connected, we must have $V=\emptyset$ and $W = \overline{\Omega}$. Moreover, $c$ is positive and $u(x)=\displaystyle\frac{L_0u(x)}{c(x)}$ is continuous and positive on $\overline{\Omega}$.
\end{dem}

\begin{lem} \label{5.1}
If $g(x)> \lambda_1$ in $\Omega$, then $L_0 u=g(x)u$ does not admit a positive solution.
\end{lem}

\begin{dem}
Suppose that $L_0 u=g(x)u$ admit a positive solution $u$ and consider $\varphi_1$ a positive eigenfunction associated to $\lambda_1$. Then,
\[
\lambda_1 \int_{\Omega}u\varphi_1 dx = \int_{\Omega}uL_0\varphi_1 dx = \int_{\Omega}\varphi_1L_0u dx = \int_{\Omega}g(x)u\varphi_1 dx
\]
that is, \[ \int_{\Omega}(g(x)-\lambda_1)u\varphi_1 dx =0,\]
since $u\varphi_1>0$ we have $g(x)-\lambda_1 =0$, which is absurd.

\end{dem}

Here, we would like to point out that in \cite{Garcia-Rossi,Garcia-Rossi 2}, García-Melián and Rossi also have considered an nonlocal  eigenvalue problem of the type 
\begin{equation}\label{eq06}
\left\{
\begin{array}{ll}
\displaystyle\int_{\R^N}K(x-y)u(y)dy- u = -\lambda u(x) & \quad \mbox{in} \ \Omega, \\
u= 0, &\quad \mbox{in} \ \R^N \backslash \Omega 
\end{array}
\right. \end{equation}
with a kernel $K \in C^1(\R^N)$, $K>0$ in $B_1$ (the unit ball), $K=0$ in $\R^N \backslash B_1$, $K(-z)=K(z)$ and $\displaystyle\int_{B_1}K(x)dx=1$. They proved that the problem (\ref{eq06}) admits a unique principal eigenvalue $\lambda_1$, that is, an eigenvalue with an associated positive eigenfunction $\phi_1 \in \C$, it is simple and verifies $0<\lambda_1(\Omega)<1$. In this way, the first eigenfunction $\phi_1$ of the problem (\ref{eq06}) is strictly positive in $\Omega$ (with a positive continuous extension to $\overline{\Omega}$) and vanishes outside $\Omega$. Therefore, a discontinuity may occur in $\partial \Omega$ and the boundary value is not taken in the usual sense, for more details see \cite[Chapter 2]{Rossi}. From these comments, we see that Proposition \ref{3} continues the study made in the above papers for another class of nonlocal problems.

\section{Framework}
 
	In whole this section, we are assuming that $K$ is a nonnegative continuous function that verifies $(K_1)$ and $(K_2)$. Moreover, for each $w \in \C$, we set the function $\Phi_{w}: \overline{\Omega} \longrightarrow \R$ by
\[
\Phi_{w}(x)= \int_{\Omega} Q(x,y)|w(y)|^{p}dy,
\]
where $p>0$ and $Q$ is a continuous function satisfying $(Q_2)$.

Since $Q$ and $w$ are bounded, we have that $\Phi_{w}$ is well defined. Furthermore, the ensuing properties will be useful later on:
\begin{itemize}
\item [$(\Phi_{1})$] $t^{p}\Phi_{w} = \Phi_{tw}$, for all $w \in \C$;
\item [$(\Phi_{2})$] $\|\Phi_{w}\|_{\infty} \leq \|Q\|_{\infty} \|w\|^{p}_{\infty}|\Omega|$, for all $w \in \C$;
\item [$(\Phi_{3})$]  $\|\Phi_{w}-\Phi_{w}\|_{\infty} \leq \|Q\|_{\infty} \||w|^{p}-|v|^{p}\|_{\infty}|\Omega|$, for all $w \in \C$;
\item [$(\Phi_{4})$] $\Phi:\C \longrightarrow \C$, given by $\Phi(w) = \Phi_{w}$, is uniformly continuous in $\C$.
\end{itemize}

Hereafter, we intend to prove the existence of positive solution for $(P)$ by using the Global Bifurcation Theorem. Having this in mind, it is very important to observe that if $(\lambda,u)$ is a solution of $(P)$, then from Lemma \ref{5}, $\lambda > \Phi_u (x)$ for all $x \in \overline{\Omega}$ (we will see that it is a necessary condition to obtain positive solution) and so
$$
L_0u=(\lambda - \Phi_u (x))u \Longleftrightarrow u=\dfrac{L_0u}{\lambda - \Phi_u (x)} \Longleftrightarrow u = \lambda^{-1}L_0u + \dfrac{\Phi_u (x)L_0u}{\lambda(\lambda - \Phi_u (x))}.
$$
For $\gamma = \lambda^{-1}$, we have 
$$
u = \gamma L_0u + \frac{\gamma^2 \Phi_u (x)L_0u}{1-\gamma \Phi_u (x)}, 
$$
or equivalently 
$$
u = \gamma L_0u + G(\gamma,u),
$$
where $G(\gamma,u)= \displaystyle\frac{\gamma^2 \Phi_u (x)L_0u}{1-\gamma \Phi_u (x)}.$

Furthermore, for each $0<a<b$,
\[
\lim_{v \rightarrow 0} \frac{G(\gamma, v)}{\|v\|_\infty} = 0, \quad \mbox{uniformly in} \ \gamma \in [a,b]. \eqno{(\cal{G})}
\]

Next, we recall the definition of compact operator when the domain is not a closed set. This type of operator applies an important role in our approach. 

\begin{defi}
Let  $A$ be an open set  in $(0,+\infty) \times \C$. An nonlinear operator  $G: {A} \longrightarrow \C$ is said to be compact if $G$ is continuous, and for each $B \subset A$ such that $B$ is bounded and $dist(B, \partial A)$ is positive, then $G(B)$ is relatively compact in $\C$.
\end{defi}


\begin{ob} \label{6}
The operator $G$ is very well defined in 
$$
A = \{(\gamma, v) \in (0,+\infty) \times \C; \gamma\|\Phi_v\|_\infty <1\}.
$$ 
Moreover, $A$ is an open set which contains $(\lambda_{1}^{-1}, 0).$ It is easy to see that $G$ is compact in each $\overline{U}_{\Lambda,\rho,M}$, where $U_{\Lambda,\rho,M} =\{(\gamma,v) \in (0,+\infty) \times \C; \|v\|_\infty <M, \Lambda^{-1}<\gamma \ \mbox{and} \ 1-\gamma\|\Phi_v\|_\infty >\rho  \}$ and $U_{\Lambda,\rho,M} \subset A$.
\end{ob}

\subsection{Proof of Theorems \ref{Teo 1} and \ref{Teo 2}}

Using the above notations, we see that $(\lambda,u)$ solves $(P)$ if, and only if,
$$L_0u + \Phi_u (x)u = \lambda u, 
$$ or equivalently, $u =F(\gamma,u):=\gamma L_0u + G(\gamma,u)$, where $\gamma = \lambda^{-1}$.

In the sequel, we will apply a Global Bifurcation Theorem found in \cite[Theorema 29.1]{Deimling}, which improves a well known Global Bifurcation Theorem found in \cite{Rabinowitz}.  

\begin{teo}\textbf{(Global bifurcation)} \label{Bifurcation}
	Let $X$ be a Banach space, $U \subset \R \times X$ a neighbourhood of $(\gamma_0,0)$, $G:\overline{U} \longrightarrow X$ completely continuous and $G(\gamma,u)=o(\|u\|_X)$ as $u \rightarrow 0$, uniformly in $\gamma$, in compacts of $\R$. Let $T \in L(X)$ be compact and $\gamma_0$ a characteristic value of odd algebraic multiplicity, $F(\gamma,u)=u-\gamma T + G(\gamma,u)$ and
	$$
	\Sigma=\{(\gamma,u)\in U; F(\gamma,u)=0,u\neq0\}.
	$$
Then the component $\mathcal{C}=\mathcal{C}_{\gamma}$ of $\overline{\Sigma}$, containing $(\gamma_0,0)$, has at least one of the following properties:
	
	(i) $\mathcal{C} \cap \partial U \neq \emptyset$
	
	(ii) $\mathcal{C}$ contains an odd number of trivial zeros $(\gamma_i,0) \neq (\gamma_0,0)$, where $\gamma_i$ is a characteristic value of $T$ of odd algebraic multiplicity.
\end{teo}

By the previous section, we know that there is a first positive eigenfunction $\varphi_{1}$ associated to $\lambda_{1}$. Moreover, $\lambda_{1}$ is an eigenvalue of $L_{0}$ with  multiplicity equal to $1$. From global bifurcation theorem, there exists a closed connected component $\mathcal{C}=\mathcal{C}_{\lambda^{-1}_1}$ of solutions for $(P)$ that satisfies $(i)$ or $(ii)$. We claim that $(ii)$ does not occur. In order to show this claim, we need of the lemma below

\begin{lem} \label{7}
	There exists $\delta>0$ such that, if $(\gamma,u)\in\mathcal{C}$ with $|\gamma-\lambda_1 ^{-1}|+\|u\|_{\infty}<\delta$ and $u\neq0$, then $u$ has defined sign, that is,
$$
u(x)>0, \quad \forall x \in \overline{\Omega} \quad \mbox{or} \quad u(x)<0, \quad \forall x \in \overline{\Omega}.
$$
\end{lem}

\begin{dem}
It is enough to prove that for any two sequences $(u_n) \subset \C$ and $\gamma_{n}\to \lambda_1^{-1}$ with 
	$$
	u_n\neq0,\quad \|u_n\|_{\infty}\to 0 \quad \mbox{and } \quad u_{n}=F(\gamma_n,u_n)=\gamma_n L_0u_n + G(\gamma_n,u_n),
$$
$u_n$ has defined signal for $n$ large enough. 

Setting $w_n=u_n/\|u_n\|_{\infty}$, we have that $(w_{n}) \subset \C$ and
\[
	w_n=\gamma_nL_0(w_n)+\frac{G(\gamma_n,u_n)}{\|u_n\|_{\infty}}=\gamma_nL_0(w_n)+o_{n}(1),
\] 
where we have used $(\cal{G})$ in the last equality.
From compactness of operator $L_0$, we can assume that $(L_0(w_n))$ is convergent for some subsequence. Then, 
$$
w_n\to w \ \mbox{in} \ \C,
$$ 
for some $w \in \C$ with $\|w\|_\infty =1$. Thereby,  
$$
w=\lambda_1^{-1}L_0(w)
$$
or equivalently,
$$
L_{0}w =\lambda_{1} w \ \mbox{in} \ \Omega.
$$ 

Thereby, $w\neq 0$ is an eigenfunction associated with $\lambda_1$, and by Proposition \ref{3} and Corollary \ref{2},  
$$
w(x)>0, \quad \forall x\in \overline{\Omega} \quad \mbox{ or } \quad w(x)<0, \quad  \forall x\in \overline{\Omega}.
$$
In the sequel, without loss of generality we assume that $w$ is positive in $\overline{\Omega}$. As $w$ is the uniform limit of $(w_{n})$ in $C(\overline{\Omega})$, we must have $w_{n}>0$ for all $x \in \overline{\Omega}$ and $n$ large enough. As $u_n$ and $w_n$ has the same signal, $u_n$ is also positive, which is the desired conclusion. 
\end{dem}

\begin{ob} 
It is easy to check that if $(\gamma,u)\in\Sigma$ if, and only if, $(\gamma,-u)$ is also in $\Sigma$. Thus, by considering the sets 
$$
\mathcal{C}^{+}=\{(\gamma,u)\in\mathcal{C}:u(x)>0,\ \forall x\in \overline{\Omega} \}\cup\{(\lambda_{1}^{-1},0)\}
$$
and
$$
\mathcal{C}^{-}=\{(\gamma,u)\in\mathcal{C}:u(x)<0,\ \forall x\in \overline{\Omega} \}\cup\{(\lambda_{1}^{-1},0)\},
$$
we have  
\begin{equation}\label{eq010}
\mathcal{C}=\mathcal{C}^{+}\cup\mathcal{C}^{-}. 
\end{equation}
Moreover, $\mathcal{C}^{-}=\{(\gamma,u)\in\mathcal{C}:(\gamma,-u)\in\mathcal{C}^{+}\}$ and $\mathcal{C}^{+}\cap\mathcal{C}^{-}=\{(\lambda_{1}^{-1},0)\}$.
\end{ob}

Indeed, in what follows, we fix
$$
\mathcal{C}^{\pm}=\{(\gamma,u)\in\mathcal{C}:u^{\pm}\not=0 \}
$$ that is, $\mathcal{C}^{\pm}$  is the subset of $\mathcal{C}$ of the functions that change signal. Since 
$$
\mathcal{C}=\mathcal{C}^+\cup\mathcal{C}^-\cup\mathcal{C}^{\pm},
$$
we deduce that to prove (\ref{eq010}), it is enough to show that $\overline{\mathcal{C}^{\pm}} = \emptyset$. Supposing by contradiction that $\overline{\mathcal{C}^{\pm}} \not= \emptyset$, as  $\mathcal{C}$ is a connected set in $(0,+\infty) \times \C$ and $\mathcal{C}^+\cup\mathcal{C}^-$ is closed nonempty set with $(\mathcal{C}^+\cup\mathcal{C}^-)\cap\mathcal{C}^{\pm}=\emptyset$, we must have 
$$
\left( \mathcal{C}^{+} \cup \mathcal{C}^{-} \right) \cap \overline{\mathcal{C}^{\pm}} \not= \emptyset.
$$ 
Therefore, there is a solution $(\gamma, u) \in \mathcal{C}$  and sequences $(\gamma_n, u_n) \subset \mathcal{C}^{+} \cup \mathcal{C}^{-} $ and $(s_n,w_n) \subset \mathcal{C}^{\pm}$ such that 
$$
\gamma_n, s_n \to \gamma \quad \mbox{in} \quad \mathbb{R} , \quad u_n \to u \quad \mbox{in} \quad \C \quad \mbox{and} \quad w_n \to u \quad \mbox{in} \quad \C. 
$$ 
Consequently $u \geq 0$ in $\overline{\Omega}$ or $u \leq 0$ in $\overline{\Omega}$ and $u \not =0$. Without loss of generality, suppose that  $u \geq 0$ in $\overline{\Omega}$. Since $(\gamma, u)$ verifies $u =\gamma L_0u+ G(\gamma,u)$, $L_0u>0$ and $G(\gamma,u) \geq 0$, it follows that $u > 0$ in $\overline{\Omega}$. Hence, $w_n$ is positive for $n$ large enough, obtaining a contradiction. Thereby,  $\overline{\mathcal{C}^{\pm}} = \emptyset$, finishing the proof of (\ref{eq010}).  

\begin{ob}
Lemma \ref{5.1} shows that the connected component that leaving $(\lambda_1,0)$, has no accumulation points of the form $(\lambda,0)$ with $\lambda > \lambda_1$.

Indeed, if $u>0$ and $\|u\|_\infty$ is small enough that $\lambda - \Phi_u (x)>\lambda_1$, from Lemma \ref{5.1}, $(\lambda,u)$ can not belong to this component.
\end{ob}

Now, consider $U \subset A$ as in Remark \ref{6} that is, $U := U_{\Lambda,\rho,M}$. Then,

\begin{lem} \label{8}
$\mathcal{C}^{+} \cap \partial U \neq \emptyset$.
\end{lem}

\begin{dem}
	Suppose by contradiction that $\mathcal{C}^{+} \cap \partial U = \emptyset$. Then, from global bifurcation theorem, there exists $(\hat{\gamma},0)\in\mathcal{C}$ , where $\hat{\gamma}\neq\lambda_{1}^{-1}$ and $\hat{\gamma}$ is a characteristic value of $L_0$ with odd algebraic multiplicity.
Hence, there exists $(u_n) \subset \C$ and $\gamma_{n}\rightarrow\hat{\gamma}$, such that
$$
u_{n}\neq0,\quad \|u_{n}\|_{\infty}\rightarrow0\quad\mbox{and}\quad u_{n}=F(\gamma_{n},u_{n}).
$$
Thus, $L_0u_n = (\lambda_n -\Phi_{u_n}(x))u_n$ and $\Phi_{u_n}(x) \rightarrow 0$ in $\Omega$, where $\lambda_n = \gamma^{-1}_n$. Then from Lemma \ref{5.1}, we have that $$\lambda_n - \Phi_{u_n}(x) > \lambda_1 \quad \mbox{for all} \ n \ \mbox{sufficiently large},$$ that is, $(\hat{\gamma},0)$ can not belong to this component, wich is an absurd.
 
\end{dem}

The next result establishes more some properties of the positive solutions of $(P)$. 

\begin{lem} \label{9}
If $(\gamma,u)$ is a solution of $u=F(\gamma,u):=  \gamma L_0u + G(\gamma,u)$ with $u>0$ in $\Omega$, then we have $$\gamma \Phi_u (x) <1 \quad \mbox{for all} \ x \in \overline{\Omega},$$ that is, $(\gamma, u) \in A=\{(\gamma, v) \in (0,+\infty) \times \C; \gamma\|\Phi_v\|_\infty <1\}$. This also states that $\mathcal{C}^{+} \subset A$.
\end{lem}

\begin{dem}
Note that, $u=F(\gamma,u):=  \gamma L_0u + G(\gamma,u)$ is equivalent to $L_0u + \Phi_u (x)=\lambda u$, where $\gamma = \lambda^{-1}$ and as $u>0$ in $\Omega$, then 
$$
0<L_0u = (\lambda-\Phi_u (x))u$$ this implies $$\lambda-\Phi_u (x)>0, \quad \mbox{for all}\ x \in \overline{\Omega}.
$$ 
Therefore, $\lambda^{-1} \Phi_u (x) <1$ for all $x \in \overline{\Omega}$ and so, 
$$
\gamma  \Phi_u (x) <1, \mbox{for all}\ x \in \overline{\Omega},
$$
that is $(\gamma,u) \in A$.
\end{dem}

Since $\overline{\Omega}$ is a compact, we can cover it with a finite number of balls centered at some points of $\overline{\Omega}$ and radius $r>0$, that is, there are $x_1,...,x_m \in \Omega$ and $m \in \N$ such that 
$$
\overline{\Omega} = \bigcup_{j=1}^{m}(B_{r/2} (x_j) \cap \overline{\Omega}),
$$ 
where $r>0$ was given in assumption $(Q_2)$. The integer $m$ will appear in the next lemma.
  
\begin{lem} \label{10}
Let $\Lambda>0$ and suppose that $(\lambda,u)$ is such that $L_0u + \Phi_u (x)u=\lambda u$ for some $\lambda \in (0, \Lambda]$. Then $\|u\|_p \leq \displaystyle\left(\frac{m \lambda}{\sigma}\right)^{\frac{1}{p}}$, that is, $u$ is uniformly bounded in $L^{p}(\Omega)$. 
\end{lem}

\begin{dem}
By the above coverage consider $E_j = \Omega \cap B_r (x_j)$, thus $$\sigma \int_{E_j}|u(y)|^pdy \leq \int_{E_j}Q(x_j,y)|u(y)|^pdy \leq \Phi_u (x_j),$$ from Lemma \ref{9}, $$\sigma \int_{E_j}|u(y)|^pdy \leq \lambda.$$ Then, $$\sigma \int_{\Omega}|u(y)|^pdy \leq m\lambda$$ that is, $\|u\|_p \leq \left(\displaystyle\dfrac{m \lambda}{\sigma}\right)^\frac{1}{p}$.
\end{dem}

As a consequence of this last proof, we have:

\begin{cor}\label{cor1}
	Under condition $(Q''_2)$. $\Phi_v(x)\geq \sigma ||v||_p^p$, for all $v\in X$, $x\in \overline\Omega$.
\end{cor}

\noindent {\bf Proof of Theorem \ref{Teo 1}}

If $u$ is a positive solution of $L_0u+\Phi_u (x)u = \lambda u$ for some $\lambda$, then $\lambda - \Phi_u (x) \geq \lambda_1$ for all $x \in \Omega$. Indeed, let $x_* \in \Omega$ be such that $\Phi_u (x_*) \leq \Phi_u (x)$ for all $x \in \Omega$. From Lemma \ref{4} we have $\|\lambda - \Phi_u\|_\infty > \lambda_1$, hence $$\lambda - \Phi_u (x_*)> \lambda_1,$$ or equivalently, $\lambda - \lambda_1 > \Phi_u (x_*)$. Moreover, from Corollary \ref{cor1}  
\begin{equation}\label{eq011}
\lambda - \lambda_1 > \Phi_u (x_*) \geq \sigma \|u\|^{p}_p. \end{equation}
On the other hand, 
\begin{equation}\label{eq012}
\lambda - \Phi_u (x) = \lambda - \Phi_u (x) + \Phi_u (x_*) - \Phi_u (x_*)$$ or $$\lambda - \Phi_u (x) > \lambda_1 -|\Phi_u (x)-\Phi_u (x_*)|. 
\end{equation}
Furthermore, from (\ref{eq011}),  
\begin{equation}\label{eq013}
|\Phi_u (x)-\Phi_u (x_*)| \leq \int_{\Omega}|Q(x,y)-Q(x_*,y)||u(y)|^pdy \leq [Q]\|u\|^{p}_{p} \leq [Q] \frac{(\lambda-\lambda_1)}{\sigma}. 
\end{equation}
Note that, if $\lambda_1 < \lambda<\lambda_1 + \dfrac{\lambda_1 \sigma}{[Q]}-\epsilon$, for a small fixed $\epsilon>0$, we have $0<\lambda -\lambda_1< \dfrac{\lambda_1 \sigma}{[Q]}-\epsilon$. Thereby, from (\ref{eq012}) and (\ref{eq013}), 
$$
\lambda- \Phi_u (x)> \lambda_1 - [Q]\dfrac{(\lambda-\lambda_1)}{\sigma}$$ and then, $$\lambda-\Phi_u(x) > \dfrac{[Q] \epsilon}{\sigma} > 0 \mbox { for all }x\in \Omega.
$$ 
Therefore, with study above we obtain that for each $\epsilon>0$, there exists $\rho >0$ ($\rho=[Q]\varepsilon\sigma^{-1}$) such that $\lambda - \|\Phi_u\|_\infty \geq \rho$, for all $\lambda \in (\lambda_1, \lambda_1 + \displaystyle\frac{\lambda_1 \sigma}{[Q]}-\epsilon]$. Moreover, 
$$
|u(x)| \leq \dfrac{|L_0u(x)|}{\lambda - \Phi_u (x)} \leq\dfrac{|L_0u(x)|}{\rho} \leq \dfrac{\|K\|_{p'}}{\rho}\|u\|_p, \quad \mbox{for all} \ x \in \overline{\Omega},
$$ 
where we have used the Hölder Inequality in the last inequality and $\displaystyle\frac{1}{p'} + \displaystyle\frac{1}{p} =1$. From Lemma \ref{10}, $\|u\|_p$ is bounded, then there exists $M>0$ such that $\|u\|_\infty \leq M$. 

Gathering all the above informations, for fixed $\epsilon>0$ and considering $\Lambda = \lambda_1 + \dfrac{\lambda_1 \sigma}{[Q]}-\epsilon$ we find $\rho>0$ and $M>0$ such that, for $U=U_{\Lambda,\frac{\rho}{2},2M}$ we have $\mathcal{C}^{+} \cap \partial U \neq \emptyset$, thus we have $(\lambda,u) \in \mathcal{C}^{+}$ and one of the following conditions occur: $\lambda -\|\Phi_u\|_\infty = \dfrac{\rho}{2}$ or $\|u\|_\infty = 2M$ or $\lambda = \Lambda$. But we have seen that, under the above conditions,  $\lambda -\|\Phi_u\|_\infty >\rho$ and $\|u\|_\infty \leq M$, that is, we should have
$\lambda =\Lambda$ and the connected component $\mathcal{C}^{+}$ crosses the hyperplane $\{\lambda\} \times \C$, for all $\lambda \in (\lambda_1,\Lambda]$. Completing the proof of the Theorem \ref{Teo 1}.

\qed

Now we will study the bifurcation for all $\lambda > \lambda_1$. In this way, as in the last proof, we want to find a positive number $\rho >0$ such that, whatever the number $\Lambda > \lambda_1$, if $u$ is a  positive solution of $L_0u+\Phi_u (x)u = \lambda u$ for some $\lambda \in (\lambda_1, \Lambda]$, then $\lambda- \|\Phi_u\|_\infty \geq \rho$. To obtain this number, we need more information on $Q$ and $p$.

\begin{lem} \label{11}
Suppose $p>0$ and $(Q_2)-(Q_3)$ hold. Then there exist $\rho >0$ and $M>0$ such that $\lambda -\|\Phi_u\|_\infty \geq \rho$ and $\|u\|_\infty \leq M$ for all $(\lambda,u)$ such that $L_0u + \Phi_u (x)u=\lambda u$, $u >0$ and $\lambda \in (\lambda_1, \Lambda]$.
\end{lem}

\begin{dem}
To begin with, we will prove first the existence of $\rho$, after we show the existence of $M$. If there is no $\rho$, then we can find a sequence $(\lambda_n,u_n)$ such that $L_0u_n +  \Phi_{u_n} (x)u_n=\lambda_n u_n$, $u_n>0$, $\lambda_n \rightarrow \lambda \in (\lambda_1,\Lambda]$ and $\|\Phi_{u_n}\|_\infty \rightarrow \lambda$. In the sequel we divide into two cases our study, namely $p >1$ and $p \in (0,1)$. \\
\noindent {\bf Case 1: ${ p>1}$.} By Lemma \ref{10},  $(u_n)$ is a bounded sequence in $L^p(\Omega)$,  and  as $p>1$, there is a some subsequence of $(u_n)$, still denoted by itself, such that  $u_n \rightharpoonup u \quad \mbox{in} \ L^p(\Omega)$. As $(L_0u_n)$ and $(\Phi_{u_n})$ are uniformly convergent in $\C$, we assume that $L_0u_n \to w$ and $\Phi_{u_n} \to v$ in $\C$ respectively. As $L_0$ is a linear and compact operator, we have 
$$
L_0u_n(x) = \int_{\Omega}K(x,y)u_n(y)dy \rightarrow \int_{\Omega}K(x,y)u(y)dy = L_0u(x) \quad \mbox{in}\ \Omega,
$$ 
and so, $L_0u_n \to L_0u$  in $\C$. Next we are going to show that $\Phi_{u_n} \to \Phi_u$ in $\C$, however as  $\Phi$ is not linear the above argument does not work well, and we need to use others arguments. From the limit $\Phi_{u_n} \to v$ in $\C$, we know that $\Phi_{u_n}(x) \rightarrow v(x), \ \forall x \in  \Omega.$ Now, as $\lambda_n-\Phi_{u_n}(x) > 0$, we have $\lambda-v(x) \geq 0$. Passing to the weak limit in the $L^p(\Omega)$ sense in $L_0u_n + \Phi_{u_n}(x)u_n = \lambda_{n}u_n$, we obtain 
$$
L_0u = (\lambda-v(x))u \quad \mbox{a.e. in}\ \Omega.
$$ 
In the sequel, we will consider the cases $u \equiv 0$ and $u\neq 0$,  and in the both cases, we will arrive to a contradiction. This proves the existence of $\rho$. 

\noindent {\bf The case $ u\equiv 0$:}  In this case, $u_n \rightharpoonup 0$ in $L^p(\Omega)$ and as $u_n>0$, we have that $u_n \rightarrow 0$ in $L^1(\Omega)$ and $L_0u_n \rightarrow 0$ in $\C$. This yields   $u_n \rightarrow 0$ in $L^p(\Omega)$. In fact, by $(Q_3)$, 
$$
\lambda_n> \Phi_{u_n}(x_0) \geq \Phi_{u_n}(x)+a(x)\|u_n\|_{p}^p \quad \mbox{for} \ x \in \Omega
$$ 
therefore, 
\begin{equation}\label{1star}
\lambda_n - \Phi_{u_n}(x) > a(x)\|u_n\|^{p}_p, \quad \mbox{for}\ x \in \Omega.
\end{equation}
From this, 
$$
\int_{\Omega} u_n(y)^p dy = \int_{\Omega}\left[\frac{L_0u_n(y)}{\lambda_n - \Phi_{u_n}(y)} \right]^pdy \leq \frac{\|L_0u_n\|^{p}_\infty}{(\|u_n\|^{p}_p)^p} \int_{\Omega}\dfrac{1}{a(y)^p}dy
$$ 
which implies
$$
\left(\int_{\Omega} u_n(y)^p dy\right)^{p+1} \leq \|L_0u_n\|^{p}_\infty\int_{\Omega}\dfrac{1}{a(y)^p}dy < \infty.
$$
Passing to the limit, and using the fact that $\|L_0u_n\|_\infty \to 0$, we find $\|u_n\|_p \rightarrow 0$. Therefore, $u_n \rightarrow 0$ in $\C$, which contradicts the limit $\|\Phi_{u_n}\|_\infty \rightarrow \lambda>0$.

\noindent {\bf The case $ u\neq 0$:} From Lemma \ref{5}, $u>0$ and $\lambda - v(x) >0$ in $\overline{\Omega}$. On the other hand, arguing as above, $u_n \rightarrow u$ in $L^{p}(\Omega)$. Hence, $\Phi_{u_n} \rightarrow \Phi_u \ \mbox{in}\ \C$ , $\|\Phi_{u_n}\|_\infty \rightarrow \|\Phi_u\|_\infty$ and $ \lambda - \Phi_u (x) >0$  for all $ \in \overline{\Omega}$. Then, $\lambda - \|\Phi_u\|_\infty >0,$ which contradicts $\|\Phi_{u_n}\|_\infty \rightarrow \lambda$.  

Now, we are going to prove the existence of $M$. Note that 
$$
|u(x)| \leq \dfrac{|L_0u(x)|}{\lambda - \Phi_u (x)} \leq\dfrac{|L_0u(x)|}{\rho} \leq \dfrac{\|K\|_q}{\rho}\|u\|_p, \quad \mbox{for all} \ x \in \overline{\Omega},
$$ 
where we have used the Hölder inequality in the last inequality with $\displaystyle\frac{1}{q} + \displaystyle\frac{1}{p} =1$. From Lemma \ref{10}, $\|u\|_p$ is bounded, then there is a $M>0$ such that $\|u\|_\infty \leq M$.

\noindent {\bf Case 2: $ p \in (0,1]$.} As in the first case, we can assume that $u_n^{p} \not\to 0$ in $L^{1}(\Omega)$, otherwise we will get $\Phi_{u_n} \to 0$ in $C(\overline{\Omega})$, which contradicts the limit $\|\Phi_{u_n}\|_\infty \rightarrow \lambda>0$. In the sequel, as $(u_n^{p})$ is bounded in $L^{1}(\Omega)$, for some subsequence, we can assume that $u_n^{p} \rightharpoonup \mu$ in $\mathcal{M}(\Omega)$ for some $\mu \in \mathcal{M}(\Omega)$, where $\mathcal{M}(\Omega)$ denotes the space of positive finite measure on $\Omega$. Thereby,  
$$
\int_{\Omega}\phi u_n^p \, dx \to \int_{\Omega}\phi d \mu, \quad \forall \phi \in C(\overline{\Omega}),
$$
and so,
$$
\Phi_{u_n}(x)=\int_{\Omega}Q(x,y) u_n^p \, dx \to \int_{\Omega}Q(x,y) d \mu=v(x), \quad \forall x \in \overline{\Omega}.
$$
As $\mu \in \mathcal{M}(\Omega)$, a  simple computation gives $v \in C(\overline{\Omega})$ and $v(x) \geq 0$ for all  $x \in \overline{\Omega}$. Using the fact that 
\begin{equation}\label{2star}
\lambda_n-\Phi_{u_n}(x) \geq a(x)\int_{\Omega}u_n^{p}\,dx,
\end{equation}
by taking the limit of $n \to +\infty$, we get  
$$
\lambda-v(x) \geq a(x)W, \quad \forall x \in \overline{\Omega}
$$
where $W=\int_{\Omega}d \mu>0$. Here we know that $W>0$, because we are supposing that $u_n^{p} \not\to 0$ in $L^{1}(\Omega)$. Since $a$ is a nonnegative function and $a^{-1} \in L^{1}(\Omega)$, we have that set $\mathcal{O}=\{x \in {\Omega} \,:\, a(x)=0\}$ has null measure, thus 
$$
\lambda -v(x)>0, \quad \mbox{a.e. in} \quad \overline{\Omega}.
$$
\begin{claim} The sequence $(u_n)$ is bounded in $L^{1}(\Omega)$.
\end{claim}

Indeed, assume by contradiction that $\|u_n\|_{1} \to +\infty$ and set $w_n=\dfrac{u_n}{\|u_n\|_{1}}$. Using the fact that $(\lambda_n,u_n)$ is a solution of $(P)$, we get
$$
L_0w_n+\Phi_{u_n}w_n=\lambda_nw_n
$$ 
and so,
$$
w_n= \frac{L_0w_n}{\lambda_n-\Phi_{u_n}}.
$$ 
As $(w_n)$ is bounded in $L^{1}(\Omega)$, for some subsequence, we have that $L_0w_n \to w_*$ in  $C(\overline{\Omega})$, consequently 
$$
w_n(x) \to \frac{w_*(x)}{\lambda-v(x)}=w(x) \quad \mbox{a.e. in} \quad \overline{\Omega}.
$$
From definition of $w$, we see that $w \in C(\Omega \setminus \mathcal{O})$ and $w(x) \geq 0$ a.e. in $\overline{\Omega}$. Moreover, we also have 
$$
w_n(x) \leq \frac{2\|L_0w_n\|_\infty}{a(x)W} \quad \mbox{a.e. in} \quad \overline{\Omega}.
$$ 
Since $a^{-1} \in L^{1}(\Omega)$, the above informations ensure that
$$
w_n \to w \quad \mbox{in} \quad L^{1}(\Omega),
$$
then $\|w\|_1=1$. On the other hand, we also have that
$$
L_0\left( \frac{u_n}{\|u_n\|_{1}^{p+1}}\right)+\Phi_{w_n}(x)w_n=\lambda_n \frac{u_n}{\|u_n\|_1^{p+1}},
$$ 
then $\Phi_{w_n}w_n \to 0$ in $L^{1}(\Omega)$, and so, by Fatou's Lemma
$$
\int_{\Omega}\int_{\Omega} Q(x,y)w^{p}(y)w^{2}(x)\,dxdy=0.
$$
From $(Q_2)$ we get $w=0$, which is an absurd. This proves that $(u_n)$ is bounded in $L^{1}(\Omega)$. Arguing as above, replacing $w_n$ by $u_n$, we can prove that $u_n \to u$ in $L^{1}(\Omega)$, and the lemma follows by repeating the same arguments explored in Case $1$.

\end{dem}

\noindent {\bf Proof of Theorem \ref{Teo 2}}

Gathering all the above informations, for all $\Lambda > \lambda_1$ we find $\rho>0$ and $M>0$ such that, for $U=U_{\Lambda,\frac{\rho}{2},2M}$ we have $\mathcal{C}^{+} \cap \partial U \neq \emptyset$, thus we have $(\lambda,u) \in \mathcal{C}^{+}$ and one of the following conditions occur: $\lambda -\|\Phi_u\|_\infty = \dfrac{\rho}{2}$ or $\|u\|_\infty = 2M$. But we have seen that, under the above conditions,  $\lambda -\|\Phi_u\|_\infty >\rho$ and $\|u\|_\infty \leq M$, that is, the connected component $\mathcal{C}^{+}$ crosses the hyperplane $\{\Lambda\} \times \C$, for all $\Lambda>\lambda_1$.

To conclude the proof, we will show the nonexistence of solution for $\lambda \leq \lambda_{1}$. In fact, suppose that $(\lambda,u)$ satisfies $u\geq 0$, $\lambda >0$ and $
L_{0}u+\Phi_{u}u=\lambda u$. Then, for all $v \in L^{2}(\Omega)$, 
\[
\langle L_{0}u,v \rangle+\langle \Phi_{u}u,v\rangle = \lambda \langle u,v \rangle.
\]
Taking $v=\varphi_{1}$, the eigenfunction associated with $\lambda_{1}$, we derive
\[
\langle L_{0}u,\varphi_{1} \rangle+\langle \Phi_{u}u,\varphi_{1}\rangle = \lambda \langle u,\varphi_{1} \rangle.
\]
As $L_{0}$ is symmetric in $L^{2}(\Omega)$, it follows that
\[
\lambda_{1} \langle \varphi_{1},u \rangle+\langle \Phi_{u}u,\varphi_{1}\rangle = \lambda \langle u,\varphi_{1} \rangle.
\]
Using the fact $\langle \Phi_{u}u, \varphi_{1} \rangle>0$, we have
\[
\lambda_{1}  \langle \varphi_{1},u \rangle< \lambda  \langle \varphi_{1},u \rangle
\]
i.e.,
\[
(\lambda_{1}-\lambda)\int_{\Omega}\varphi_{1}udx <0,
\]
showing that $\lambda_{1}<\lambda$.

\section{Proof of Theorem \ref{Teo 3}}

In this section, we will fix $\lambda>\lambda_1$ and show that Problem $(P)$ has a positive solution that will be the uniform limit of solutions given by Theorem \ref{Teo 2}.

First, consider $0< \epsilon \leq \epsilon_0 = \dfrac{N}{2p}$ and define 
\[
Q_\epsilon (x,y)=Q(x,y)(2-a_\epsilon (x)) 
\]
 where 
$$
a_\epsilon (x) = \left\{
\begin{array}{cll}
|x-x_0|^\epsilon, \quad \mbox{if} \ |x-x_0| \leq 1 \\
1, \quad \mbox{if} \ |x-x_0| \geq 1 .
\end{array}
\right.
$$
Note that, $2Q(x,y) \geq Q_\epsilon (x,y) \geq Q(x,y)\geq 0$ for all $x,y \in \Omega$ and $Q_\epsilon (x,y) \geq \sigma$ when $|x-y| \leq r$, and so, $Q_\epsilon$ verifies $(Q_2)$. Moreover, 
\[
Q_\epsilon (x_0,y) - Q_\epsilon (x,y) = 2Q(x_0,y) - Q(x,y)(2-a_\epsilon (x))\\
\geq 2Q(x,y) - Q(x,y)(2-a_\epsilon (x))
\]
that is, 
\begin{equation}\label{eq015}
Q_\epsilon (x_0,y) - Q_\epsilon (x,y) \geq Q(x,y)a_\epsilon (x)\geq\frac 12 Q_\epsilon (x,y)a_\epsilon(x), \quad \mbox{for all} \ x,y \in \Omega.
\end{equation} 
 
Related to $a_\epsilon$ we have: $a_\epsilon(x)\leq1$ for all $x\in \Omega$.

We will consider the following family of auxiliary problems:
$$
L_0u + \Phi^{\epsilon}_u (x)u = \lambda u, \quad \mbox{in} \ \Omega. \eqno{(P_\epsilon)}
$$

\begin{lem} \label{12}
Suppose that $(Q_2)$, $(Q_4)$ hold and fix $\epsilon >0$. Then,  if $(\lambda, u)$ is a positive solution of $(P_\epsilon)$ with $\lambda>\lambda_1$ and $u>0$, we have 
\begin{equation}\label{star}
\lambda - \Phi^{\epsilon}_u (x) \geq \theta a_\epsilon (x), \quad \mbox{for all} \ x \in \Omega,
\end{equation}
 where $\theta=\min\{\lambda_1,\lambda-\lambda_1\}$.
\end{lem}

\begin{dem}
Note that, $(Q_2), (Q_4)$ and (\ref{eq015}) imply in $\Phi^{\epsilon}_u (x_0) \geq \Phi^{\epsilon}_u (x)$, for all $x \in \Omega$ and, 
\begin{equation}\label{eq017}
\lambda- \Phi^{\epsilon}_u (x) \geq \Phi^{\epsilon}_u (x_0)-\Phi^{\epsilon}_u (x) \geq a_\epsilon (x) \int_\Omega Q(x,y) u(y)^p dy \geq \frac12a_\epsilon (x) \Phi^{\epsilon}_u (x). \end{equation} 
Thus, if $\Phi^{\epsilon}_u (x) \leq \lambda_1$, we have $\lambda - \Phi^{\epsilon}_u (x) \geq \lambda - \lambda_1\geq (\lambda-\lambda_1)a_\epsilon(x)$. 
On the other hand, if $\Phi^{\epsilon}_u (x) > \lambda_1$, from (\ref{eq017}) we have $\lambda - \Phi^{\epsilon}_u (x) \geq \lambda_1 a_\epsilon (x)$. Anyway, since $\theta=\min\{\lambda_1,\lambda-\lambda_1\}$ we have that
\[
\lambda - \Phi^{\epsilon}_u (x) \geq \theta a_\epsilon (x), \quad \mbox{for all} \ x \in \Omega.
\]

\end{dem}

As a consequence of Lemma \ref{12}, we see that the same conclusion of Lemma \ref{11} for the auxiliary problem $(P_\epsilon)$. 

\begin{lem} \label{121}
	Suppose that $0<\epsilon\leq\epsilon_0$ is fixed, $p>0$,  $(Q_2)$ and $(Q_4)$ hold. Then there exist $\rho >0$ and $M>0$ such that $\lambda -\|\Phi^\epsilon_u\|_\infty \geq \rho$ and $\|u\|_\infty \leq M$ for all $(\lambda,u)$ satisfying $(P_\epsilon)$, with $u >0$ and $\lambda \in (\lambda_1, \Lambda]$.
\end{lem}

\begin{dem} We proceed exactly as in the proof of Lemma \ref{11}. Here we replace (\ref{1star}) by (\ref{star}), in the case $p>1$, and we replace (\ref{2star}) by (\ref{star}) for the case $p\leq1$.
	
\end{dem}

Using Lemma \ref{121} and following all the steps in the proof of Theorem \ref{Teo 2}, for any fixed $\epsilon\in (0,\epsilon_o]$,  we have  a connected component $\mathcal C_\epsilon^+$, associated to the bifurcation equation $(P_\epsilon)$, which crosses the hyperplans $\{\lambda\}\times \C$ for all $\lambda>\lambda_1$.

\begin{ob}\label{exist}
We resume this last remark as follows: for any $\lambda>\lambda_1$ and any $\epsilon\in (0,\epsilon_0]$ we have a positive $u\in \C$ satisfying $(P_\epsilon)$.
\end{ob}

Now we are ready to prove Theorem \ref{Teo 3}.

\begin{dem} (\textbf{Theorem \ref{Teo 3})} Fix $\lambda>\lambda_1$ again and consider the functions $g_n:\overline{\Omega} \longrightarrow \R$ defined by $$g_n (x)=\int_{\Omega} Q_\frac{1}{n}(x,y)|u_n(y)|^pdy,$$ where $u_n$ is given by Remark \ref{exist}, with $\epsilon=\frac1n$, which verifies $L_0u_n + g_n (x)u_n = \lambda u_n$. The proof consists in proving that the problem $(P)$ has a solution which is a limit of a subsequence of $u_n$ when $n$ goes to infinity.

We know that $\|u_n\|_p$ is bounded and $g_n$ is bounded in $L^\infty (\Omega)$. Moreover, $g_n$ is uniformly convergent in compact parts of $\overline{\Omega} \setminus \{x_0\}$. In fact, 
\[
|g_n (x)-g_m (x)| \leq 2\|Q\|_\infty \|u_n\|^{p}_p \big| |x-x_0|^{\frac{1}{n}} - |x-x_0|^{\frac{1}{m}} \big|, \quad \mbox{if} \ |x-x_0| \leq 1
\] and
\[
|g_n (x)-g_m (x)|=0, \quad \mbox{if} \ |x-x_0| > 1.
\] 

In the sequel we divide into two cases our study, namely $p >1$ and $p \in (0,1)$. \\
\noindent {\bf Case 1: $p>1$.} By Lemma \ref{10},  $(u_n)$ is a bounded sequence in $L^p(\Omega)$,  and  as $p>1$, there is a some subsequence of $(u_n)$, still denoted by itself, such that  $u_n \rightharpoonup u \quad \mbox{in} \ L^p(\Omega)$. As in the proof of Theorem \ref{Teo 2}, $L_0u_n$ converges to  $L_0u$ uniformly in $\overline \Omega$.  Denote by $v$ the uniform limit of $g_n $ in the compact parts of $\overline\Omega\setminus \{x_0\}$. 
Let us show that, $u_n$ converges to $u$ in $L^p(\Omega)$. From Lemma \ref{12}  we have
\[
\lambda - g_n(x) \geq \theta a_\frac{1}{n}(x) \geq \theta R^{\frac{1}{n}}, \quad \mbox{for} \ x \in \overline{\Omega}\setminus B_R(x_0), \ (R<1)
\]
that is, $\lambda - g_n(x)$ converges uniformly to $\lambda - v(x)$ in $\overline{\Omega}\setminus B_R(x_0)$ which implies $\lambda - v(x) \geq \theta$ in $\overline{\Omega}\setminus B_R(x_0)$, hence 
\begin{equation}\label{eq0016}
u_n(x)=\dfrac{L_0 u_n(x)}{\lambda - g_n(x)} \rightarrow \dfrac{L_0u(x)}{\lambda - v(x)} \quad \mbox{uniformly in} \ \overline{\Omega}\setminus B_R(x_0).
\end{equation} Moreover, $v(x)<\lambda$ for $x \neq x_0$. Now, in the neighborhood of $x_0$, we have
\begin{equation}\label{eq0017}
\int_{B_R(x_0)} u_n(y)^p dy = \int_{B_R(x_0)}\left[\frac{L_0u_n(y)}{\lambda - g_n(y)} \right]^pdy \leq \frac{\|L_0u_n\|^{p}_\infty}{\theta^p} \int_{B_R(x_0)}\dfrac{1}{a_\frac{1}{n}(y)^p}dy \leq \dfrac{\|L_0u_n\|^{p}_\infty \omega_N R^{N-\frac{p}{n}}}{\theta^p (N-\frac{p}{n})}.
\end{equation}
The $L^p$-convergence of $(u_n)$ follows from (\ref{eq0016}) and (\ref{eq0017}).

Now, we claim that $g_n$ converges to $\Phi_u$ in $\overline{\Omega}$. In fact, since $(u_n)$ converges in $L^p(\Omega)$, passing to a subsequence if necessary, $(u_n)$ is dominated by a function $h \in L^p(\Omega)$ and $$Q_\frac{1}{n}(x,y) u_n(y)^p \leq 2\|Q\|_\infty h(y)^p.$$ The assertion follows from the dominated convergence.
As $g_n(x)<\lambda$, we have $\Phi_u(x) \leq \lambda$. Passing to the weak limit in the $L^p(\Omega)$ sense in $L_0u_n + g_n(x)u_n = \lambda u_n$, we obtain 
$$
L_0u = (\lambda-\Phi_u(x))u \quad \mbox{a.e. in}\ \Omega.
$$ 
In the sequel, we will consider the cases $u \equiv 0$ and $u\neq 0$.

\noindent {\bf The case $ u\equiv 0$:} We can not have $u \equiv 0$, because $g_n(x)\leq 2 \|Q\|_\infty \|u_n\|^{p}_p$ and thus $g_n$ converges uniformly to $0$ in $\overline{\Omega}$, which contradicts the Lemma \ref{5.1} because we would have $\lambda -g_n(x)>\lambda_1$ in $\overline{\Omega}$ for $n$ large enough.

\noindent {\bf The case $u\neq 0$:} From Lemma \ref{5}, $u>0$ and $\lambda - v(x) >0$ in $\overline{\Omega}$. On the other hand, arguing as above, $u_n \rightarrow u$ in $L^{p}(\Omega)$. Hence, $g_n \rightarrow \Phi_u \quad \mbox{in}\ \C$ , $\|g_n\|_\infty \rightarrow \|\Phi_u\|_\infty$ and $ \lambda - \Phi_u (x) >0$  for all $ \in \overline{\Omega}$. Then $u$ is a positive solution sought.  

\noindent {\bf Case 2: $ p \in (0,1]$:} As in the first case, we can assume that $u_n^{p} \not\to 0$ in $L^{1}(\Omega)$. In the sequel, as $(u_n^{p})$ is bounded in $L^{1}(\Omega)$, for some subsequence, we can assume that $u_n^{p} \rightharpoonup \mu$ in $\mathcal{M}(\Omega)$ for some $\mu \in \mathcal{M}(\Omega)$, where $\mathcal{M}(\Omega)$ denotes the space of positive finite measure on $\Omega$. Thereby,  
$$
\int_{\Omega}\phi u_n^p \, dy \to \int_{\Omega}\phi d \mu, \quad \forall \phi \in C(\overline{\Omega}),
$$
and so,
$$
g_n(x)=\int_{\Omega}Q_\frac{1}{n}(x,y) u_n^p \, dy \to \int_{\Omega}Q(x,y) d \mu=v(x), \quad \forall x \in \overline{\Omega}.
$$
As $\mu \in \mathcal{M}(\Omega)$, a  simple computation gives $v \in C(\overline{\Omega})$ and $v(x) \geq 0$ for all  $x \in \overline{\Omega}$. Using the fact that 
$$
\lambda_n-g_n(x) \geq \theta a_\frac{1}{n}(x), \quad \mbox{for} \ x \neq x_0
$$
by taking the limit of $n \to +\infty$, we get  
$$
\lambda-v(x) \geq \theta >0, \quad \forall x \in \overline{\Omega}\setminus \{x_0\}
$$
Hence, $\lambda-v(x) \geq \theta >0,\ a.e.$ in $\overline{\Omega}$.

\begin{claim} The sequence $(u_n)$ is bounded in $L^{1}(\Omega)$.
	
\end{claim}
Indeed, assume by contradiction that $\|u_n\|_{1} \to +\infty$ and set $w_n=\dfrac{u_n}{\|u_n\|_{1}}$. Using the fact that $(\lambda,u_n)$ is a solution of $(P)$, we get
$$
L_0w_n+g_n(x)w_n=\lambda w_n
$$ 
and so,
$$
w_n= \frac{L_0w_n}{\lambda_n-g_n(x)}.
$$ 
As $(w_n)$ is bounded in $L^{1}(\Omega)$, for some subsequence, we have that $L_0w_n \to w_*$ in  $C(\overline{\Omega})$, consequently 
$$
w_n(x) \to \frac{w_*(x)}{\lambda-v(x)}=w(x) \quad \mbox{a.e. in} \quad \overline{\Omega}.
$$
From definition of $w$, we see that $w \in C(\Omega)$ and $w(x) \geq 0$ a.e. in $\overline{\Omega}$. Moreover, we also have 
$$
w_n(x) \leq \frac{2\|L_0w_n\|_\infty}{\theta} \quad \mbox{a.e. in} \quad \overline{\Omega}.
$$ 
The above informations ensure that
$$
w_n \to w \quad \mbox{in} \quad L^{1}(\Omega),
$$
then $\|w\|_1=1$. On the other hand, we also have that
$$
L_0\left( \frac{u_n}{\|u_n\|_{1}^{p+1}}\right)+\Phi^{\frac{1}{n}}_{w_n}(x)w_n=\lambda \frac{u_n}{\|u_n\|_1^{p+1}},
$$
where $$\Phi^{\frac{1}{n}}_{w_n} (x) = \int_{\Omega} Q_{\frac{1}{n}} (x,y) |w_n(y)|^p dy.$$
Then $\Phi^{\frac{1}{n}}_{w_n}w_n \to 0$ in $L^{1}(\Omega)$, and so, by Fatou's Lemma
$$
\int_{\Omega}\int_{\Omega} Q(x,y)w^{p}(y)w^2(x)\,dxdy=0,
$$
from $(Q_2)$ we get $w=0$, which is absurd. This proves that $(u_n)$ is bounded in $L^{1}(\Omega)$. Arguing as above, replacing $w_n$ by $u_n$, we can prove that $u_n \to u$ in $L^{1}(\Omega)$, and the result follows by repeating the same arguments explored in Case $1$. 

\end{dem}

\textbf{Final comments:} We expect that the bifurcation results in this work remain valid for weaker conditions on $Q$, without assumptions $(Q_3)$ and $(Q_4)$ for instance.

To finish this work, we would like to remark that the same result of this section holds under the condition $(Q'_4)$:
\begin{itemize}
	\item[$(Q'_4)$]  There are a decomposition $\Omega=\displaystyle{\bigcup_{j=1}^m E_j}$, $x_1\in E_1$, $x_2\in E_2$,... $x_m\in E_m$  such that $Q(x_j,y) \geq Q(x,y)$ for all $x\in E_j$, $y \in \Omega$.
\end{itemize}

It is easy to check if we consider $Q_\epsilon$ replacing $a_\epsilon$ by 
$$
a_\epsilon (x) = \left\{
\begin{array}{cll}
|x-x_1|^\epsilon....|x-x_m|^\epsilon, \quad \mbox{if} \ |x-x_0| \leq 1 \\
1, \quad \mbox{if} \ |x-x_0| \geq 1 ,
\end{array}
\right.
$$
then the proof works with the following adjusts. The uniform convergence of $g_n$ is in the compact parts of $\overline\Omega\setminus\{x_1,...,x_m\}$. In the case 1, $p>1$, $u_n$ converges uniformly in $$\overline\Omega\setminus\bigcup_{j=1}^m B_R(x_j), \mbox { for small } R>0 $$
and for any $j=1,2,...,m$,
\[
\int_{B_R(x_j)} u_n(y)^p dy = \int_{B_R(x_j)}\left[\frac{L_0u_n(y)}{\lambda - g_n(y)} \right]^pdy \leq \frac{\|L_0u_n\|^{p}_\infty}{\theta^p} \int_{B_R(x_j)}\dfrac{1}{a_\frac{1}{n}(y)^p}dy \leq \dfrac{\|L_0u_n\|^{p}_\infty \omega_N R^{N-\frac{p}{n}}}{\theta^p (N-\frac{p}{n})}.
\]
The proof follows as in the proof of Theorem \ref{Teo 3}.

\end{document}